\documentstyle[fleqn,12pt]{article}
\begin{document}\pagenumbering{arabic}\setcounter{page}{1}
\pagestyle{plain}\baselineskip=16pt

\thispagestyle{empty}
\rightline{MSUMB 97-01, January 1997} 
\vspace{1.4cm}

\begin{center}
{\Large\bf Differential calculus on the $h$-superplane}
\end{center}

\vspace{0.5cm}
\noindent
Salih Celik \\
{\footnotesize Department of Mathematics, Mimar Sinan University, 80690 Besiktas, Istanbul, TURKEY. }

\noindent
Sultan A. Celik\\
{\footnotesize Department of Mathematics, Yildiz Technical University, Sisli, Istanbul, TURKEY. }

\noindent
Metin Arik \\
{\footnotesize Department of Physics, Bogazici University, Bebek, Istanbul, TURKEY. }

\vspace{1.5cm}
{\bf Abstract}

A non-commutative differential calculus on the $h$-superplane is 
presented via a contraction of the $q$-superplane. An R-matrix 
which satisfies both ungraded and graded Yang-Baxter equations is 
obtained and a new deformation of the $(1+1)$ dimensional classical 
phase space (the super-Heisenberg algebra) is introduced. 

\vfill\eject
\noindent
{\bf I. INTRODUCTION}

A possible approach to quantum groups is obtained by deforming the 
coordinates of a linear space to be non-commuting objects.$^1$ 
In this scheme the quantum group structure appears if one considers 
linear transformations which preserve the algebraic properties of the 
algebra of coordinates. A natural and physically interesting question 
is whether one can define differentials and derivatives corresponding 
to these non-commuting variables. A general answer to this question 
was given by Connes.$^2$ Who considered the differential algebra for 
non-commutative algebras. Differential geometry of Lie groups and 
supergroups plays an important role in the mathematical modeling of 
physical theories. Since a (graded) Hopf algebra or quantum (super) 
group can be regarded as a generalization of the notion of a (super) 
group, it is tempting to also generalize the corresponding notions of 
differential geometry.$^3$ 

Recently Wess and Zumino$^4$ have shown that a consistent quantum 
deformation of the differential calculus is satisfied by an 
R-matrix which can be any solution of the quantum Yang-Baxter 
equation. The non-commutative plane has been studied in the course 
of the last few years$^5$ in that the geometry underlying the above 
mentioned plane has a deep connection with the Yang-Baxter equation 
which is important in two dimensional exactly soluble statistical 
models. The quantum (super) plane was generalized to the supersymmetric 
quantum (super) plane$^{6,7}$ and it was shown that the 
quantum superplane is related to the graded Yang-Baxter equation. 

The $h$-deformation of the supergroup GL$(1|1)$ is given by Dabrowski 
and Parashar$^8$, via a contraction of GL$_q(1|1)$. They have also 
introduced a differenetial calculus on the $h$-superplane which is related 
to GL$_h(1|1)$ via the Wess-Zumino formulae. 

This paper considers an alternative approach where instead of adopting 
the Wess-Zumino calculus from the start, the R-matrix is obtained using 
the consistency conditions. This leads to a consistent exterior derivative. 
We start by introducing a non-commutative differential calculus 
on the $h$-superplane via a contraction of the $q$-superplane. We 
define derivatives and differentials on the $h$-superplane of 
non-commuting coordinates and give their commutation rules. We 
note the role of the graded Yang-Baxter equation and the connection 
to the quantum supergroup. We give a new deformation of the $(1+1)$ 
dimensional classical phase space. Finally we show that the $q$-deformed 
super-oscillator algebra satisfies the undeformed (classical) 
super-oscillator algebra when objects are transformed into new objects 
such that they are singular for certain values of the deformation 
parameter. 

\noindent
{\bf II. A DIFFERENTIAL CALCULUS ON $h$-SUPERPLANE} 

In this work we denote $q$-deformed objects by primed quantities. 
Unprimed quantities represent transformed coordinates. As usual, 
we assume that even (bosonic) objects commute with everything and 
odd (grassmann) objects anticommute among themselves. 

\noindent
{\bf A. Quantum $h$-superplane }

We begin by considering the quantum superplane which is defined by 
Manin.$^9$ The commutation relation between the even coordinate $x'$ 
and the odd (grassmann) coordinate $\theta'$ of the quantum superplane 
is in the form 
$$ x' \theta' - q \theta' x' = 0, \eqno(1) $$
where $q$ is a complex deformation parameter. 

We now introduce new coordinates $x$ and $\theta$, in terms of $x'$ 
and $\theta'$ as 
$$x = x', \quad \theta = \theta' - {h\over{q - 1}} x' \eqno(2)$$
as in Ref. 10. This transformation is singular in the 
$q \longrightarrow 1$ limit. Using relation (1), it is easy to verify that 
$$x \theta = q \theta x + h x^2, \eqno(3)$$
where the new deformation parameter $h$ commutes with the coordinate $x$. 
Also, since the grassmann coordinate $\theta'$ satisfies 
$$\theta'^2 = 0$$
one obtains 
$$\theta^2 = - h \theta x, \eqno(4)$$
where $h$ anti-commutes with $\theta$ and 
$$h^2 = 0, \eqno(5)$$
that is, the new deformation parameter $h$ is a grassmann number.$^8$ 
Taking the $q \longrightarrow 1$ limit we obtain the following relations 
which define the $h$-superplane 
$$x \theta = \theta x + h x^2, \quad 
\theta^2 = - h \theta x. \eqno(6)$$

\noindent
{\bf B. Relations of Coordinates and Differentials} 

To establish a non-commutative differential calculus on the quantum 
$h$-superplane, we assume that the commutation relations between the 
coordinates and their differentials are in the following form 
$$x' {\sf d}x' = A {\sf d}x' x' $$
$$x' {\sf d}\theta' = F_{11} {\sf d}\theta' x' + F_{12} {\sf d}x' \theta', 
  \eqno(7)$$
$$\theta' {\sf d}x' = F_{21} {\sf d}x' \theta' + F_{22} {\sf d}\theta' x', $$
$$\theta' {\sf d}\theta' = B {\sf d}\theta' \theta'. $$

We demand an exterior differential {\sf d} obeying the condition: 
$${\sf d}^2 = 0, \eqno(8)$$
and the graded Leibniz rule 
$${\sf d}(f g) = ({\sf d} f) g + (- 1)^{\hat{f}} f ({\sf d} g), \eqno(9)$$
where $\hat{f}$ is the grassmann degree of $f$ (recall that {\sf d} should be 
odd), that is, $\hat{f} = 0$ for even variables and $\hat{f} = 1$ 
for odd variables. Considering the differential of a function and 
differentiating (2) we have 
$$ {\sf d}x = {\sf d}x', \quad {\sf d}\theta = {\sf d}\theta' + 
  {h\over {q - 1}} {\sf d}x'. \eqno(10)$$
Substituting (2) and (10) into (7) one has 
$$x~ {\sf d}x = A {\sf d}x~ x, $$
$$x~ {\sf d}\theta = F_{11} {\sf d}\theta~ x + F_{12} {\sf d}x~ \theta + 
  {h\over {q-1}}\left( 
  A - F_{11} - F_{12}\right) {\sf d}x~ x, \eqno(11)$$
$$\theta~ {\sf d}x = F_{21} {\sf d}x~ \theta + F_{22} {\sf d}\theta~ x - 
  {h\over {q-1}}\left(
  A + F_{21} + F_{22}\right) {\sf d}x~ x, $$
$$ \theta~ {\sf d}\theta = B {\sf d}\theta~ \theta -  {h\over {q-1}}\left[
  \left(B + F_{12} + F_{21}\right) {\sf d}x~ \theta -  
  \left(B - F_{11} - F_{22}\right) {\sf d}\theta~ x \right]. $$

These relations are slightly different from the results of Dabrowski and 
Parashar. $^8$ The reason for this difference is that in ref. 8, the 
commutation relations among the matrix elements of a matrix belonging to 
$GL_h(1\vert 1)$ were obtained via the use of commutation relations of the 
dual {\it exterior} superplane. On the other hand we use the commutation 
relations of the dual superplane as in ref. 3. The commutation relations 
among the matrix elements are the same in both approaches. However the 
commutation relations involving the differentials turn out not to be the 
same. 

We know that the quantum supermatrices of the quantum $GL_q(1|1)$ 
supergroup can be defined as linear transformations of the variables $x'$ 
and $\theta'$ which preserve the commutation relation (1) and their 
duals$^9$, that is, the quantum supergroup $GL_q(1|1)$ acts as a 
linear transformation on the quantum superplane which preserves (1) and 
the relations 
$$\varphi'^2 = 0, \quad \varphi' y' - q^{-1} y' \varphi' = 0. \eqno(12)$$
In extending this property of covariance under the coaction of 
$GL_q(1|1)$, from the superplane to its calculus, it will be assumed 
that the deformed group structure implies and is implied by invariance 
of the intermediary relations (7) under linear transformations of the 
quantum superplane. In the present work, this will be applied to the 
$h$-deformed superplane.  

The coefficients $A, B$ and $F_{ij}$ given in (11) can be related to $q$ 
by the consistency of calculus. Thus we apply the exterior derivative {\sf d} 
to the relation (3). From the consistency condition 
$${\sf d}(x \theta - q \theta x - h x^2) = 0 \eqno(13)$$
we find 
$$F_{11} = q (1 - F_{22}), \quad F_{12} = - (1 + q F_{21}).\eqno(14\mbox{a})$$
Applying the exterior derivative {\sf d} on the second and third 
relations of (11) and using the definitions (8), (9) one finds 
$$F_{12} = q F_{11} - 1, \quad F_{21} = q (F_{22} - 1). \eqno(14\mbox{b})$$
Similarly from the last relation of (11) we get 
$$F_{12} + F_{21} = q (F_{11} + F_{22}) - (1 + q) B. \eqno(14\mbox{c})$$
On the other hand, using the relation (4) we obtain 
\begin{eqnarray*}
\theta~ {\sf d}\theta 
& = & {\sf d}\theta~ \theta + h F_{21} {\sf d}x~ \theta - 
      h (1 - F_{22}) {\sf d}\theta~ x \\
& = & {\sf d}\theta~ \theta - {h\over q} \left[F_{11} {\sf d}\theta~ x + 
      (1 + F_{12}) {\sf d}x~ \theta \right].
\hspace*{4.2cm}(15)
\end{eqnarray*}
Thus we have 
$$B = 1, \quad F_{12} + F_{21} + 1 = (1 - q)F_{21}, $$
$$1 - F_{11} - F_{22} = (1 - q)(1 - F_{22}). \eqno(14\mbox{d})$$

To find the commutation relation between differentials, say ${\sf d}x$ and 
${\sf d}\theta$, we apply the exterior derivative {\sf d} on the first two 
relations of (11) and use the nilpotency of {\sf d} [eq. (8)]. Then it 
is easy to see that 
$${\sf d}x ~{\sf d}\theta = {1\over q} {\sf d}\theta ~{\sf d}x, \quad 
  \left({\sf d}x\right)^2 = 0.\eqno(16)$$
Solving the system (14) we now get 
$$A ~~\mbox{undetermined}, \quad F_{11} = q,$$
$$F_{12} = q^2 - 1, \quad F_{21} = - q, \quad F_{22} = 0. \eqno(17)$$
We choose $A$ equal to $q^2$, since this leads to the standard R-matrix 
[eq. (38)] in the $h \longrightarrow 0$ limit. We are thus led to the 
following deformed relations containing $q$ and $h$ 
$$ x ~{\sf d}x = q^2 {\sf d}x~ x, $$
$$ x~ {\sf d}\theta = q {\sf d}\theta~ x - h {\sf d}x~ x + 
   (q^2 - 1) {\sf d}x~ \theta, \eqno(18)$$
$$\theta~ {\sf d}x = - q {\sf d}x~ \theta - q h {\sf d}x~ x, $$
$$ \theta~ {\sf d}\theta = {\sf d}\theta~ \theta - 
   h \left(q {\sf d}x~ \theta + {\sf d}\theta~ x \right).$$
Note that although in the $q \longrightarrow 1$ limit the transformations 
(2) and (10) are ill behaved, the resulting commutation relations are well 
defined. We shall not use the limit prosess, yet. 

\noindent
{\bf C. Relations of Derivatives and Variables} 

First, we wish to consruct the curl of any one-form $w(x',\theta')$. 
To this end let us denote the partial derivatives with 
respect to $x'$ and $\theta'$ as 
$$\partial_{x'} = {\partial\over {\partial x'}}, \quad 
  \partial_{\theta'} = {\partial\over {\partial \theta'}}, \eqno(19)$$
respectively, and introduce the super-gradiend operator in 
vector notation 
$$\nabla = \left(\partial_{x'}, \partial_{\theta'}\right). \eqno(20)$$
We define the vectors $X' = (x', \theta')$ and 
${\sf d}X' = ({\sf d}x', {\sf d}\theta')$. 
Then we can write the differential {\sf d} as 
$${\sf d} = {\sf d}x' \partial_{x'} + {\sf d}\theta' \partial_{\theta'} = 
  {\sf d}X'.\nabla \eqno(21)$$
where the dot denotes the inner product. 
If we write $w(x',\theta')$ in the basis ${\sf d}X'$ as 
$$w(x',\theta') = {\sf d}x' w_1(x',\theta') + {\sf d}\theta' w_2(x',\theta') 
  \eqno(22)$$
where $w_1$ and $w_2$ are smooth functions of the variables then we get 
$${\sf d}w(x',\theta') = {\sf d}x' {\sf d}\theta' 
  \left[q \partial_{x'} w_2(x',\theta') - 
  \partial_{\theta'} w_1(x',\theta') \right]. \eqno(23)$$
Thus the curl of the one-form 
$w(x',\theta') = \left(w_1(x',\theta'), w_2(x',\theta')\right)$ 
is given by 
$$\nabla \times w(x',\theta') = q \partial_{x'} w_2(x',\theta') - 
  \partial_{\theta'} w_1(x',\theta'). \eqno(24)$$
Now we can find the commutation rules of the derivatives, once we obtain 
the derivatives $\partial_{x}$ and $\partial_{\theta}$. For this if we 
demand the chain rule on the expressions (2) we find 
$$\partial_{x} = \partial_{x'} + {h\over {q - 1}} \partial_{\theta'}, \quad 
  \partial_{\theta} = \partial_{\theta'}. \eqno(25)$$
It is easy to see that in the case of (25) the differential {\sf d} given in 
(21) preserves its form 
$${\sf d} = {\sf d}x \partial_{x} + {\sf d}\theta \partial_{\theta}. 
  \eqno(26)$$

If we now put $w(x,\theta) = {\sf d} f(x,\theta)$ then we conclude that 
$\partial_{x}$ and $\partial_{\theta}$ generate a non-commutative algebra 
with the commutation relations 
$$\partial_\theta \partial_x = q \partial_x \partial_\theta, \quad 
  \partial_{\theta}^2 = 0. \eqno(27)$$
From (26) 
$${\sf d} f(x,\theta) = {\sf d}x \partial_{x}f + 
  {\sf d}\theta \partial_{\theta}f, \eqno(28)$$
so that replacing $f$ with $xf$ and $\theta f$ we arrive at the following 
commutation relations between derivatives and variables 
$$ \partial_x x = 1 + A x \partial_x + F_{12} \theta \partial_\theta - 
  {h\over {q - 1}} (A - F_{11} - F_{12}) x \partial_\theta, $$
$$ \partial_x \theta = - F_{21} \theta \partial_x - {h\over {q - 1}} 
  \left[(A + F_{21} + F_{22}) x \partial_x + (1 + F_{12} + F_{21}) 
  \theta \partial_\theta \right], $$
$$ \partial_\theta x = F_{11} x \partial_\theta, \eqno(29)$$
$$\partial_\theta \theta = 1 - \theta \partial_\theta - F_{22} x \partial_x 
   - {h\over {q -1}}(1 - F_{11} - F_{22}) x \partial_\theta. $$
These commutation relations are well defined in the limit 
$q \longrightarrow 1$. This can be checked using (17). 

The covariant differential structure which is obtained so far is a 
concrete example of non-commutative differential geometry. The 
complete framework of the differential calculus requires commutation 
relations of the exterior differentials with derivatives. 

\noindent
{\bf D. Relations of Differentials with Derivatives}

Finally we shall find the commutation relations between differentials 
and derivatives. We assume that they have the following form in terms 
of primed quantities 
$$ \partial_{x'} {\sf d}x' = A_{11} {\sf d}x' \partial_{x'} + 
   A_{12} {\sf d}\theta' \partial_{\theta'},$$
$$ \partial_{x'} {\sf d}\theta' = A_{21} {\sf d}\theta' \partial_{x'} + 
   A_{22} {\sf d}x' \partial_{\theta'}, \eqno(30)$$
$$ \partial_{\theta'} {\sf d}x' = B_{11} {\sf d}x' \partial_{\theta'} + 
   B_{12} {\sf d}\theta' \partial_{x'}, $$
$$\partial_{\theta'} {\sf d}\theta' = B_{21} {\sf d}\theta' \partial_{\theta'} 
  + B_{22} {\sf d}x' \partial_{x'}. $$
Using (25) and (10), these commutation rules can be written as 
$$ \partial_x {\sf d}x = A_{11} {\sf d}x \partial_x + {h\over {q - 1}} 
   (A_{11} - A_{12} + B_{11}) {\sf d}x \partial_\theta + 
   A_{12} {\sf d}\theta \partial_\theta + 
   {h\over {q - 1}} B_{12} {\sf d}\theta \partial_x, $$
$$\partial_x {\sf d}\theta = A_{21} {\sf d}\theta \partial_x + 
  {h\over {q - 1}} (A_{11} - A_{21} + B_{22}) {\sf d}x \partial_x + 
 {h\over {q - 1}} (A_{12} - A_{21} + B_{12}) {\sf d}\theta \partial_\theta\, $$
$$\partial_\theta {\sf d}\theta = B_{21} {\sf d}\theta \partial_\theta + 
      {h\over {q - 1}}(B_{22} - B_{21} - B_{11}) {\sf d}x \partial_\theta + 
     B_{22} {\sf d}x \partial_x - 
     {h\over {q - 1}} B_{12} {\sf d}\theta \partial_x, $$
$$\partial_\theta {\sf d}x = B_{11} {\sf d}x \partial_\theta + 
   B_{12} {\sf d}\theta \partial_x - {h\over {q - 1}} B_{12} 
   ({\sf d}\theta \partial_\theta + {\sf d}x \partial_x). \eqno(31)$$
Finding the coefficients $A_{ij}$ and $B_{ij}$ from these equations we get 
$$ \partial_{x} {\sf d}x = q^2 {\sf d}x \partial_{x} - 
   h {\sf d}x \partial_{\theta} + (q^2 - 1) {\sf d}\theta \partial_{\theta}, $$
$$ \partial_{x} {\sf d}\theta = q {\sf d}\theta \partial_x + 
   q h \left({\sf d}x \partial_x + {\sf d}\theta \partial_\theta\right), $$
$$ \partial_\theta {\sf d}x = - q {\sf d}x \partial_{\theta}, \eqno(32)$$
$$ \partial_\theta {\sf d}\theta = {\sf d}\theta \partial_{\theta} + 
   h {\sf d}x \partial_{\theta}. $$
Here, in order to obtain these relations we used that the exterior 
differential {\sf d} (anti-) commutes with the differentials, that is, 
$${\sf d}~ ({\sf d}x) = - ({\sf d}x)~ {\sf d}, \quad 
  {\sf d}~ ({\sf d}\theta) = ({\sf d}\theta)~ {\sf d} \eqno(33)$$
and the relation 
$$\partial_i (X^j {\sf d}X^k) = \delta^i{}_j \delta^k{}_l {\sf d}X^k 
  \eqno(34)$$
where $\partial_1 = \partial_x$, $\partial_2 = \partial_\theta$, $X^1 = x$, 
and $X^2 = \theta$. Again, the relations (32) are well defined in the limit 
$q \longrightarrow 1$. 

\noindent
{\bf E. The R-Matrix Formalism} 

We now shall obtain the R-matrix satisfying the graded Yang-Baxter 
equations 
$$R_{12} R_{13} R_{23} = R_{23} R_{13} R_{12}, \eqno(35)$$ 
$$\hat{R}_{12} \hat{R}_{23} \hat{R}_{12} = 
  \hat{R}_{23} \hat{R}_{12} \hat{R}_{23}. \eqno(36)$$ 
For this we define the commutation relations between variables and their 
differentials [see eq. (18)] in the following form 
$$X^i {\sf d}X^j = q (- 1)^{\hat{i}(\hat{j} + 1)} K^{ji}{}_{kl} 
  {\sf d}X^k X^l \eqno(37)$$
where $K \in End({\cal C}\otimes {\cal C})$. Comparing (37) with (18) we have 
$$K_{h,q}  
 = \left(\matrix{   q          &  0          & 0   & 0 \cr 
                    h          &  1          & 0   & 0 \cr
                    - q^{-1} h &  q - q^{-1} & 1   & 0 \cr 
                    0   &  - h & - q^{-1} h  & q^{-1}  \cr }\right) 
  = (K^{ij}{}_{kl}). \eqno(38)$$
If we define 
$$K_h = \lim_{q\rightarrow 1} K_{h,q}, \eqno(39)$$
$$\hat{K}_h = \lim_{q\rightarrow 1} (K_{h,q} P), \eqno(40)$$
where $P$ is the super permutation matrix, that is, 
$$P^{ij}{}_{kl} = (- 1)^{\hat{i} \hat{j}} \delta^i{}_l \delta^j{}_k,$$
we have 
$$K_h = \left(\matrix{   
    1 &  0   &   0   & 0 \cr 
    h &  1   &   0   & 0 \cr
  - h &  0   &   1   & 0 \cr 
    0 &  - h & - h   & 1  \cr }\right), \quad  
\hat{K}_h = \left(\matrix{   
  1  &  0   & 0    & 0 \cr 
  h  &  0   & 1    & 0 \cr
- h  &  1   & 0    & 0 \cr 
  0  &  - h & - h  & - 1  \cr }\right). \eqno(41)$$
Here the matrix $\hat{K}_h$ coincides with the $\hat{R}_h$ matrix 
of ref. 8. We know from ref. 8 that the matrix $\hat{K}_h$ satisfies 
equation (36) with the grading 
$$(\hat{K}_{12})^{abc}{}_{def} = \hat{K}^{ab}{}_{de} \delta^c{}_f, $$
$$(\hat{K}_{13})^{abc}{}_{def} = (-1)^{b(c + f)} 
  \hat{K}^{ac}{}_{df} \delta^b{}_e, $$
$$(\hat{K}_{23})^{abc}{}_{def} = (-1)^{a(b + c + e + f)} 
  \hat{K}^{bc}{}_{ef} \delta^a{}_d. 
  \eqno(42)$$
Also, the R-matrix 
$$R_h = P \hat{K}_h \eqno(43)$$
obeys both the ungraded and the graded Yang-Baxter equations with the 
grading again given by (42). This is due to the odd character of $h$. 
As a consequence, 
$$K_h = R_h^{-1} \eqno(44)$$
so that $K_h$ has the same properties as $R_h$. 

We note that the equation 
$$\hat{K}_h T_1 T_2 = T_1 T_2 \hat{K}_h \eqno(45)$$
is satisfied for the $h$-deformed supergroup $GL(1|1)$. Here $T$ is a 
supermatrix in $GL_h(1|1)$ and $T_1 = T \otimes I$, $T_2 = I \otimes T$. 
It is assumed that the tensor product is graded, that is, 
$$(T_1)^{ij}{}_{kl} = T^i{}_k \delta^j{}_l, \quad 
  (T_2)^{ij}{}_{kl} = (- 1)^{\hat{i}(\hat{j} + \hat{l})} 
  T^j{}_l \delta^i{}_k. $$

It is well-known that a supermatrix has in the form 
$$T = \left(\matrix{ a & \beta \cr \gamma & d \cr}\right)$$
with two even (latin letters) and two grassmann (greek letters) matrix 
elements. Equation (45) explicitly reads 
$$a \beta = \beta a, \quad a \gamma = \gamma a + h (a^2 + \gamma \beta - ad),$$
$$d \beta = \beta d, \quad d \gamma = \gamma d - h (d^2 - \gamma \beta - da),$$
$$ \beta^2 = 0, \quad  \gamma^2 = h \gamma (d - a), $$
$$ \beta \gamma = - \gamma \beta + h \beta (d - a),$$
$$ a d = d a + h \beta (a - d). \eqno(46)$$

Using the $K_h$ matrix, we now formulate to the differential calculus on 
the $h$-superplane. 
The commutation relations between variables and their differentials are 
$$X^i {\sf d}X^j = (- 1)^{\hat{i}(\hat{j} + 1)} K^{ji}{}_{kl} 
  {\sf d}X^k X^l. \eqno(47)$$
The commutation relations between variables and derivatives are 
$$\partial_j X^i = \delta^i{}_j + (- 1)^{\hat{i} \hat{j}} K^{ik}{}_{lj} 
   X^l \partial_k, \eqno(48)$$
and the relations between differentials and derivatives are 
$$\partial_j {\sf d}X^i = (-1)^{\hat{i}(\hat{j} + 1)} (K^{-1})^{ik}{}_{lj} 
   {\sf d}X^l \partial_k. \eqno(49)$$
Note that the commutation relations between variables can be expressed 
using the $\hat{K}$ matrix as 
$$X^i X^j = \hat{K}^{ij}{}_{kl} X^k X^l, \eqno(50)$$
and the relations between derivatives as 
$$\partial_i \partial_j = \hat{K}^{kl}{}_{ji} 
\partial_l \partial_k. \eqno(51)$$

\noindent
{\bf F. The Commutation Relations} 

We would like to discuss the meaning of covariance in a graded version of 
non-commutative differential calculus of Wess-Zumino. $^4$ Before proceeding, 
we define the dual quantum $h$-superplane. For this we interpret the 
differentials ${\sf d}x$ and ${\sf d}\theta$, as the coordinates of the dual 
superplane as follows [see eq. (12)] 
$${\sf d}x = \varphi, \quad {\sf d}\theta = y. \eqno(52)$$
We now formulate the differential calculus on the $h$-superplane 
as follows: 

The commutation relations of variables and their differentials are 
$$x \theta = \theta x + h x^2, \quad \theta^2 = - h \theta x, $$
$$\varphi y = y \varphi, \quad \varphi^2 = 0. \eqno(53\mbox{a})$$

Note that if we assume that the relations (53a) have to be covariant under 
the coaction 
$$\delta(x) = a \otimes x + \beta \otimes \theta, \qquad 
  \delta(\theta) = \gamma \otimes x + d \otimes \theta, $$
and that $\beta$, $\gamma$ anticommute with $\theta$, $\varphi$ and $h$ 
we get anew the relations (46). 

The commutation relations among the derivatives are 
$$\partial_x \partial_\theta = \partial_\theta \partial_x, \quad 
  \partial_\theta^2 = 0, \eqno(53\mbox{b})$$
and those between variables and derivatives are 
$$\partial_x x = 1 + x \partial_x + h x \partial_\theta, \quad 
  \partial_\theta x = x \partial_\theta, $$
$$\partial_x \theta = \theta \partial_x - 
  h (x \partial_x + \theta \partial_\theta), \quad 
  \partial_\theta \theta = 1 - \theta \partial_\theta + h x \partial_\theta. 
  \eqno(53\mbox{c})$$

The commutation relations of variables with their differentials are 
$$x~ {\sf d}x = {\sf d}x~ x, \quad 
  \theta~ {\sf d}\theta = {\sf d}\theta~ \theta - 
   h ({\sf d}x~ \theta + {\sf d}\theta~ x),$$
$$x~ {\sf d}\theta = {\sf d}\theta~ x - h {\sf d}x~ x, \quad 
  \theta~ {\sf d}x = - {\sf d}x~ \theta - h {\sf d}x~ x. \eqno(54)$$

The commutation relations between derivatives and differentials are 
$$\partial_x {\sf d}x = {\sf d}x \partial_x - h {\sf d}x \partial_\theta, 
  \quad \partial_\theta {\sf d}\theta = {\sf d}\theta \partial_\theta + 
  h {\sf d}x \partial_\theta, $$
$$\partial_x {\sf d}\theta = {\sf d}\theta \partial_x + 
  h ({\sf d}x \partial_x + {\sf d}\theta \partial_\theta), \quad 
  \partial_\theta {\sf d}x = - {\sf d}x \partial_\theta. 
\eqno(55) $$

Note that the simple relations involving the differentials $\varphi$ and $y$ 
in (53) are not obtained using the R-matrix $\hat{K}_h$. Instead the 
$q \longrightarrow 1$ limit of (16) has been used. 

The exterior differential 
$$ {\sf d} = \varphi \partial_x + y \partial_\theta \eqno(56)$$
satisfies the usual properties such that 
$$ {\sf d} x - x {\sf d} = \varphi, \quad 
   {\sf d} \theta + \theta {\sf d} = y, \eqno(57\mbox{a})$$
as expected. The relations of the exterior differential {\sf d} with 
$\partial_x$ and $\partial_\theta$ are 
$$ {\sf d} \partial_x = \partial_x {\sf d}, \quad 
   {\sf d} \partial_\theta = - \partial_\theta {\sf d} \eqno(57\mbox{b}) $$
and those with differentials $\varphi$ and $y$ 
$${\sf d} \varphi = - \varphi {\sf d}, \quad {\sf d} y = y {\sf d}. 
  \eqno(57\mbox{c}) $$
Thus the basic requirement for the exterior derivative is quite consistent: 
\begin{eqnarray*} 
{\sf d}^2 & = & {\sf d} (\varphi \partial_x + y \partial_\theta) \\
    & = & - \varphi {\sf d} \partial_x + y {\sf d} \partial_\theta \\
    & = & - (\varphi \partial_x + y \partial_\theta) {\sf d} = - {\sf d}^2 
\end{eqnarray*}
so that ${\sf d}^2$ must vanish.

As we understand from the language of Wess-Zumino $^4$, covariance here 
means that all the relations between coordinates $x$, $\theta$, 
differentials ${\sf d}x$, ${\sf d}\theta$ and derivatives 
$\partial_x$, $\partial_\theta$, etc. must preserve their form when one 
changes the coordinates by 
$$ x \longrightarrow a x + \beta \theta, $$
$$ \theta \longrightarrow \gamma x + d \theta, \eqno(58\mbox{a})$$
where the matrix $T = (T^i_j)$ is an element of the quantum supergroup 
acting on the quantum superspace. 
Therefore we change the differentials by 
$$ {\sf d}x \longrightarrow a {\sf d}x - \beta {\sf d}\theta, $$
$${\sf d}\theta \longrightarrow - \gamma {\sf d}x + d {\sf d}\theta. 
   \eqno(58\mbox{b})$$
This is consistent since the exterior differential {\sf d} anticommutes with 
the grassmann variables as mentioned before. 
Covariance can be maintained if one defines the transformation law of the 
partial derivatives as folllows 
$$\partial_x \longrightarrow (a^{-1} - a^{-1} \gamma d^{-1} \beta a^{-1}) 
  \partial_x - a^{-1} \gamma d^{-1} \partial_\theta, $$
$$\partial_\theta \longrightarrow (d^{-1} - d^{-1} \beta a^{-1} \gamma d^{-1}) 
  \partial_\theta + d^{-1} \beta a^{-1} \partial_x.  \eqno(58\mbox{c}) $$
Hence the differentials transform under the action of supertranspose of 
supertranspose of $T$, $(T^{st})^{st}$, whereas the derivatives transform 
under the inverse of the supertranspose of $T$, $(T^{st})^{-1}$. 

\noindent
{\bf III. A NEW DEFORMATION OF CLASSICAL PHASE SPACE}

We now shall obtain a new deformation of the $(1 + 1)$-dimensional 
super-Heisenberg algebra (the classical phase space). We 
denote the algebra (53) generated by coordinates $x$, $\theta$ 
and the derivatives $\partial_x$ and $\partial_\theta$ by ${\cal B}_h$. 
It is interesting to note that simply identifying 
$\partial_x$ and $\partial_\theta$ with $i p_x$ and $p_\theta$ 
is not compatible with the hermiticity of coordinates and 
momenta. In other words, the algebra ${\cal B}_h$ 
cannot be interpreted as a deformation of the $(1 + 1)$-dimensional 
super-Heisenberg algebra. In order to identify $\partial_x$ and 
$\partial_\theta$ with the momenta $i p_x$ and $p_\theta$, one must 
take care of hermiticity of coordinates and momenta. To this end, 
let us define the hermitean conjugation of the coordinates 
$x$ and $\theta$, respectively, as 
$$x^+ = x, \quad \theta^+ = \theta + 2h x. \eqno(59)$$
It is then easy to see that the deformation parameter of the algebra 
(53) becomes a pure imaginary parameter: 
$$\overline{h} = - h \eqno(60)$$
where the bar denotes complex conjugation. In this case, the 
hermitean conjugation of the derivatives $\partial_x$ and $\partial_\theta$ 
are 
$$\partial_x^+ = - \partial_x + 2h \partial_\theta, \quad 
  \partial_\theta^+ = \partial_\theta. \eqno(61)$$
Note that the definitions (59) and (61) is for the classical case are 
obtained in the $h \longrightarrow 0$ limit. 

The relations (53) are now invariant under the definitions (59)-(61). The above 
involution allows us to define the hermitean operators 
$$\hat{x} = x, \quad \hat{\theta} = \theta + h x, \eqno(62)$$
and 
$$\hat{p}_x = i (\partial_x - h \partial_\theta), \quad 
  \hat{p}_\theta = \partial_\theta. \eqno(63)$$
The final form of the $h$-deformed super-Heisenberg algebra is 
$$\hat{x} \hat{\theta} = \hat{\theta} \hat{x} + 
  h \hat{x}^2, \quad 
  \hat{\theta}^2 = - h \hat{\theta} \hat{x}, $$
$$\hat{p}_x \hat{p}_\theta = \hat{p}_\theta \hat{p}_x, \quad 
   \hat{p}_\theta^2 = 0, $$
$$\hat{p}_x \hat{x} = \hat{x} \hat{p}_x + 
   i (1 + h \hat{x} \hat{p}_\theta), \quad 
  \hat{p}_\theta \hat{x} = \hat{x} \hat{p}_\theta, \eqno(64)$$
$$\hat{p}_x \hat{\theta} = \hat{\theta} \hat{p}_x - 
  h (\hat{x} \hat{p}_x + i \hat{\theta} \hat{p}_\theta), \quad 
 \hat{p}_\theta \hat{\theta} = 
  1 - \hat{\theta} \hat{p}_\theta + 
  h \hat{x} \hat{p}_\theta.$$

\noindent
{\bf IV. A COMMENT ON SUPER-OSCILLATORS} 

We know that introducing one 'bosonic' and one 'fermionic' oscillator, $A$ 
and $B$, respectively, and making the usual identification 
$$x' ~\longleftrightarrow~ A^+, \quad \theta' ~\longleftrightarrow~ B^+,$$
$$\partial_{x'} ~\longleftrightarrow~ A, \quad 
  \partial_{\theta'} ~\longleftrightarrow~ B, \eqno(65)$$
one constructs the quantum super-oscillator algebra which is covariant under 
the quantum supergroup GL$_q(1|1)$. Under identification (25) and (2) 
give 
$$x ~\longleftrightarrow~ A^+, \quad 
  \partial_x ~\longleftrightarrow~ A + {h\over {q - 1}} B,$$
$$ \partial_\theta ~\longleftrightarrow~ B, \quad 
  \theta ~\longleftrightarrow~ B^+ - {h\over {q - 1}} A^+ \eqno(66)$$
where $q$ is a real number. Substituting (66) into (29) and using (17), 
surprisingly all $h$-dependence cancels and one obtains the usual 
$q$-deformed super-oscillator algebra$^{11}$ 
$$A A^+ = 1 + q^2 A^+ A + (q^2 - 1) B^+ B, $$
$$B B^+ = 1 - B^+ B, \quad B^2 = 0 = B^{+ 2}, $$
$$A B^+ = q B^+ A, \quad A B = q^{- 1} B A. \eqno(67)$$
In the $q \longrightarrow 1$ limit, we get undeformed 
super-oscillator algebra. 

An interesting problem is the construction of a differential calculus 
on the $(h_1,h_2)$-superplane using the methods of this paper. 
A differential calculus on this superplane using the Wess-Zumino 
formulae has been given in Ref. 12. 

\noindent
{\bf ACKNOWLEDGEMENT}

This work was supported in part by {\bf T. B. T. A. K.} the 
Turkish Scientific and Technical Research Council. 

We would like to express our deep gratitude to the referee for critical 
comments on the manuscript. 

\noindent
$^1$ V. G. Drinfeld, {\em Quantum groups}, 
    in {\it Proc. } IMS, Berkeley, (1986); \\
\hspace*{0.18cm}{
    Yu I. Manin, {\it Quantum groups and non-commutative geometry}, CRM,} \\
\hspace*{0.38cm}{Montreal University, (1988). }\\
$^2$ A. Connes, {\it Non-commutative differential geometry}, Publ. Math. 
     I.H.E.S. {\bf 62} 
     \hspace*{0.18cm} (1985). \\
$^3$ B. Schmidke, S. Vokos and B. Zumino, {\it Z. Phys. C} {\bf 48}, 
     249 (1990); \\
     \hspace*{0.18cm} E. Corrigan, B. Fairlie, P. Fletcher and R. Sasaki, 
     {\it J. Math. Phys.} {\bf 31}, 776 
     \hspace*{0.18cm} (1990). \\
$^4$ Wess, J. and Zumino, B., {\it Nucl. Phys. B } (Proc. Suppl.) 
    {\bf 18} B, 302 (1990).\\
$^5$ F. Mueller-Hoisen, {\it J. Phys. A} {\bf 25}, 1703 (1992).\\
$^6$ S. Soni, {\it J. Phys. A} {\bf 24}, L459 (1991); 
     ibid 619. \\
$^7$ W.S. Chung, {\it J. Math. Phys.} {\bf 35}, 2484 (1994). \\
$^8$ L. Dabrowski and P. Parashar, {\it Lett. Math. Phys.} {\bf 38}, 
     (1996).\\
$^9$ Yu. I. Manin, {\em Commun. Math. Phys.} {\bf 123}, 163 (1989). \\
$^{10}$ A. Aghamohammadi, M. Khorrami and A. Shariati, 
        {\it J. Phys. A} {\bf 28}, L225 
        \hspace*{0.28cm} (1995);\\
      \hspace*{0.28cm} V. Karimipour, {\it Lett. Math. Phys.} {\bf 30}, 
      87 (1994). \\
$^{11}$ M. Chaichian, P. Kulish and J. Lukierski, {\it Phys. Lett. B} 
        {\bf 262}, 43 (1991). \\ 
$^{12}$ Salih \c Celik, {\it Lett. Math. Phys.} {\bf 42}, 299 (1997). 

\end{document}